\documentclass{amsart}

\usepackage{amsmath}
\usepackage{amsthm}
\usepackage{amssymb}
\usepackage[all]{xy}

\UseComputerModernTips

\newcommand{\nat}{\ensuremath{\mathbb{N}}}

\newcommand{\ganz}{\ensuremath{\mathbb{Z}}}
\newcommand{\R}{\ensuremath{\mathbb{R}}}
\newcommand{\C}{\ensuremath{\mathbb{C}}}
\newcommand{\Rp}{\ensuremath{\mathbb{R}_{\ge 0}}}

\newcommand{\eins}{1 \hspace{-2.3pt} \mathrm{l}}

\newcommand{\tate}{\ensuremath{\mathbb{L}}}

\newcommand{\Knull}{\ensuremath{\operatorname{K}_0}}
\newcommand{\KnullG}{\ensuremath{\operatorname{K}_0^G}}
\newcommand{\KnullGstr}{\ensuremath{\operatorname{K^{\prime}}_0^G}}

\newcommand{\Knullgrstr}[1]{\ensuremath{\operatorname{K^{\prime}}_0^{#1}}}

\newcommand{\Vark}{\ensuremath{\mathit{Var}_k}}
\newcommand{\VarS}{\ensuremath{\mathit{Var}_S}}
\newcommand{\Var}[1]{\ensuremath{\mathit{Var}_{#1}}}

\newcommand{\Bl}[2]{\ensuremath{\operatorname{Bl}_{#1}#2}}

\newcommand{\cS}{\ensuremath{\mathcal{S}}}
\newcommand{\M}{\ensuremath{\mathcal{M}}}
\newcommand{\D}{\ensuremath{\mathcal{D}}}
\newcommand{\reg}{\ensuremath{\mathcal{O}}}

\newcommand{\Tau}{\ensuremath{\mathrm{T}}}
\newcommand{\aff}[1]{\ensuremath{\mathbb{A}^{#1}}}
\newcommand{\proj}[1]{\ensuremath{\mathbb{P}^{#1}}}
\newcommand{\Proj}{\ensuremath{\operatorname{\mathbb{P}}}}
\newcommand{\Spec}{\ensuremath{\operatorname{Spec}}}

\newcommand{\Ind}{\ensuremath{\operatorname{Ind}}}
\newcommand{\muedach}{\ensuremath{\hat{\mu}}}
\newcommand{\multgr}{\ensuremath{\mathbb{G}_m}}
\newcommand{\bog}{\ensuremath{\operatorname{\mathcal{L}}}}
\newcommand{\ord}{\ensuremath{\operatorname{ord}}}
\newcommand{\De}{\ensuremath{\mathbb{D}}}
\newcommand{\sd}{\ensuremath{\;|\;}}
\newcommand{\Star}{\ensuremath{\operatorname{Star}}}
\newcommand{\Stab}{\ensuremath{\operatorname{Stab}}}

\newcommand{\jac}[1]{\ensuremath{\mathcal{J}_{#1}}}
\newcommand{\raus}{\ensuremath{\backslash}}
\newcommand{\irr}{\ensuremath{\operatorname{irr}}}

\newcommand{\Pl}{\ensuremath{\mathcal{P}}}

\newtheorem{thm}{Theorem}[section]

\newtheorem{lem}[thm]{Lemma}
\newtheorem{prop}[thm]{Proposition}
\newtheorem{cor}[thm]{Corollary}

\newtheorem{listprop}[thm]{List of properties}
\theoremstyle{definition}
\newtheorem{defi}[thm]{Definition}
\theoremstyle{remark}
\newtheorem{rem}[thm]{Remark}

\newtheorem{claim}[thm]{Claim}

\numberwithin{equation}{section}

\begin{document}
\title{On motivic zeta functions and the motivic nearby fiber}
\author{Franziska Bittner}
\address{Universit\"at Essen, FB6, Mathematik,45117 Essen, Germany}
\email{franziska.bittner@uni-essen.de}
\keywords{Motivic zeta function, motivic nearby fiber}
\subjclass{14F42, 32S30}
\begin{abstract}
We collect some properties of the motivic zeta functions and the
motivic nearby fiber defined by Denef and Loeser. In particular, we
calculate the relative dual of the motivic nearby fiber. We give a
candidate for a nearby cycle morphism on the level of Grothendieck
groups of varieties using the motivic nearby fiber.
\end{abstract}
\maketitle
\section{Introduction}
Let $k$ be an algebraically closed field of characteristic zero.
One parameter Taylor series of length $n$ in a smooth variety $X$ over $k$
are called \emph{arcs of order $n$} on $X$. The set of all these arcs
are the $k$-valued points of a $k$-variety $\bog_n(X)$.

Suppose that we are given a function
$f:X\longrightarrow \aff{1}$ on a smooth connected variety $X$ of
dimension $d$.
Denef and Loeser have associated to this data the so-called \emph{motivic
zeta function} $S(f)(T)$, which is a formal power series with coefficients
in a localized equivariant Grothendieck group of varieties over the
zero locus of $f$. 
The $n$-th coefficient of this series (for $n\ge 1$) is
given as $\tate^{-nd}$ times the class of the variety of arcs $\gamma(t)$
of order $n$ on $X$ satisfying $f\circ\gamma(t)=t^n$, which carries a 
$\mu_n$-action induced by $t\mapsto\zeta t$.
Here $\tate$ denotes the class of the affine line.  
Using the transformation rule for motivic integrals, they have given 
a formula for $S(f)(T)$ in terms of an embedded
resolution of the zero locus, which shows that it is in fact a
rational function which is regular at infinity. Let us denote
$-S(f)(\infty)$, the \emph{motivic nearby fiber},  
by $\psi_f$, as it is supposed 
to be a virtual motivic incarnation of the nearby cycle sheaf.
We establish some identities for the motivic nearby fiber
which are analogues of identities known
for the nearby cycle sheaf. In particular, we calculate the relative
dual over the zero locus of $f$. It turns out to be $\tate^{1-d}\psi_f$,
which means that $\psi_f$ behaves as the class of a smooth variety, proper over
the zero locus. 
Furthermore we give a functional equation for the motivic zeta function.
Unfortunately, to be able to define e.g. a relative dual, 
we have to work in Grothendieck groups which are
coarser than the ones considered by Denef and Loeser.

For the purpose of investigating the zeta function and the nearby fiber, we
first have to generalize slightly the presentations from~\cite{euler} in the
equivariant setting. Essentially, we also allow \emph{free} actions on
the base varieties. Then we calculate the relative dual of an
affine toric variety given by a simplicial cone which is proper over
the base variety.

Using the motivic nearby fiber, we end with defining a 
nearby cycle morphism on the level of Grothendieck groups of varieties
and listing some properties of this morphism.

I am indebted to Eduard Looijenga, my thesis advisor. I thank Jan
Denef and Wim Veys for their comments and questions.

\subsection*{Conventions:}
In the sequel $k$ denotes an algebraically closed field of
characteristic zero. By a \emph{variety} over $k$ we mean a reduced
scheme of finite type over $k$.

\section{Motivic zeta functions and the motivic nearby fiber}
\label{motint}
The aim of this section is to recall the definition and properties of
motivic zeta functions and the motivic nearby fiber as they can be
found e.g.\ in \cite{DenefLoeser} and in \cite{Looijenga}.
We also use this opportunity to fix some notations. 

\subsection{Arc spaces}
We denote $\Spec(k[t]/t^{n+1})$ by $\De_n$.
An \emph{order $n$ arc} in a $k$-variety $X$ is a morphism
$\gamma:\De_n\longrightarrow X$. 
There is a $k$-scheme $\bog_n(X)$,
the \emph{space of arcs of order $n$} in $X$, whose $k$-valued points
are the arcs of order $n$ on $X$. Actually, it is the scheme
representing the functor from $k$-schemes to sets which sends $Z$ to
$X(Z\times \De_n)$. There are natural projections
$\bog_{n+1}(X)\longrightarrow \bog_n(X)$ which are $\aff{\dim X}$-bundles
in case $X$ is smooth and equidimensional.

\subsection{Grothendieck groups of varieties} \label{Grothendieck}
We denote by 
$\Knull(\Vark)$ the free abelian group on isomorphism classes of
$k$-varieties modulo the relations $[X]=[X-Y]+[Y]$, where $Y\subset
X$ is a closed subvariety. It is also called the \emph{naive
Grothendieck group of varieties over $k$} and carries a ring structure
induced by the product of varieties.
More generally, for a variety $S$ over $k$, we denote by
$\Knull(\VarS)$ the free abelian group on isomorphism classes $[X]_S$
of varieties over $S$ modulo the relations $[X]_S=[X-Y]_S+[Y]_S$, 
where $Y\subset X$ is a closed subvariety. It is naturally a
$\Knull(\Vark)$-module. We denote by $\tate$ the class of the affine
line $[\aff{1}]\in\Knull(\Vark)$. Let $\M_S$ be the localization
$\Knull(\VarS)[\tate^{-1}]$ and $\M_k$ the ring
$\Knull(\Vark)[\tate^{-1}]$. 
A morphism $f:S\longrightarrow S'$ naturally induces 
$f_!:\M_S\longrightarrow \M_{S'}$ by composition and
$f^*:\M_{S'}\longrightarrow \M_S$ by pulling back. The presentation
for $\Knull(\VarS)$ given in \cite{euler} allows to define a relative
dualization endomorphism $\D_S$ on $\M_S$ characterized by the
property $\D_S[X]_S=\tate^{-\dim X}[X]_S$ for a smooth equidimensional
variety $X$ which is proper over $S$. 

Let $S$ be a variety over $k$ with a good
action of a finite group $G$. Here 
by a \emph{good} action we mean an action such that
every orbit is contained in an affine open subvariety.
By $\KnullGstr(\Var{S})$ we denote the free group on isomorphism
classes $[X]_S$ of varieties $X\longrightarrow S$ with a good
$G$-action over the action on $S$ modulo relations for closed
$G$-subvarieties.
It has the structure of a $\KnullGstr(\Vark)$-module provided by the
product of varieties with the diagonal action.

We define $\KnullG(\Var{S})$ to be $\KnullGstr(\Var{S})$ modulo the
submodule generated as a group by expressions of the form
$[G \circlearrowright \Proj(V)]_S
-[\proj{n}\times (G\circlearrowright X)]_S$, where 
$V\longrightarrow X$ is a vector bundle of rank $n+1$ with a
$G$-action which is linear over the action on $X$.
Here $G \circlearrowright \Proj(V)$ denotes the projectivization of
this action, whereas \linebreak
$\proj{n}\times (G\circlearrowright X)$ denotes the action of $G$ on
the right factor only. For $S=\Spec k$ this is an ideal, and 
$\KnullG(\VarS)$ is a $\KnullG(\Vark)$-module. 
For an $\aff{n}$-bundle $F\longrightarrow X$ with an
affine $G$-action over the action on the base we have
$[G\circlearrowright F]_S=[\aff{n}\times(G\circlearrowright X)]_S$ 
in $\KnullG(\Var{S})$.
We obtain restriction and induction morphisms from group
homomorphisms. In particular, $\KnullG(\VarS)$ is a $\Knull(\Vark)$-module.
Denote by $\M_S^G$ 
the localization $\Knull^G(\VarS)[\tate^{-1}]$.
An equivariant morphism $f:S\longrightarrow S'$ induces 
$f_!:\M^G_S\longrightarrow \M^G_{S'}$ and
$f^*:\M^G_{S'}\longrightarrow \M^G_S$. 
In case of a trivial action of $G$ on the base variety, the
presentation in~\cite{euler} allows to define a relative dualization
endomorphism $\D_S$ on $\M^G_S$ characterized by the
property $\D_S[X]_S=\tate^{-\dim X}[X]_S$ for a smooth equidimensional
variety $X$ which is proper over $S$. 

Now let $\muedach=\varprojlim \mu_m$ (where in the projective
limit we take the natural surjections $\mu_{dm}\longrightarrow \mu_m,
\zeta\mapsto \zeta^d$). 
Following~\cite{DenefLoeser}, we call an action of $\muedach$ which
comes from a good action of a finite quotient $\mu_n$ of $\muedach$ 
a \emph{good $\muedach$-action}.

As in the case of a finite group for a base variety $S$ (with a
trivial $\muedach$-action) we can then define 
the equivariant group $\Knull^{\muedach}(\VarS)$. Actually, it is just
the direct limit
$\varinjlim\Knull^{\mu_m}(\VarS)$, where 
$\Knull^{\mu_m}(\VarS)\longrightarrow\Knull^{\mu_{dm}}(\VarS)$ is restriction
along the $d$-th power map $\mu_{dm}\longrightarrow \mu_m$.
We denote $\Knull^{\muedach}(\VarS)[\tate^{-1}]$ by $\M_S^{\muedach}$.

\begin{rem}
The equivariant groups which Denef and Loeser consider are finer. The
groups we use
allow us to define dualization, induction morphisms, 
quotient morphisms by groups acting freely on the base variety and
zeta-functions in an equivariant setting.  
\end{rem}

\subsection{Motivic zeta functions}

Let $X$ be a smooth connected variety over $k$ of dimension $d$, 
let $f:X\longrightarrow \aff{1}$ be a function. 
For a natural number $n\ge 0$ we define

$$\cS_n(f):=\tate^{-nd}[\{\gamma\in\bog_n(X)\,|\,\ord_t(f\gamma)=n\}]_X\in
\M_X.$$

\begin{rem}
Sometimes, for example in~\cite{MotivicIgusa}, there is an extra
factor $\tate^{-d}$.
\end{rem}

Following Denef and Loeser, we 
define the \emph{naive motivic zeta function of $f$} as 
$$\cS(f)(T):=\sum_{n\ge 0}\cS_n(f)T^n \in \M_X[[T]].$$

\begin{rem}
Sometimes (for example in \cite{DenefLoeser}) the constant term
$\cS_0(f)$ is omitted. We keep it, like in~\cite{MotivicIgusa} and 
\cite{Looijenga}.
\end{rem}

Denef and Loeser established 
a formula for 
$\cS(f)(T)$ in terms of an
embedded
resolution $H:Y\longrightarrow X$ of $f^{-1}(0)$. Let us assume that
$f$ is non-zero.
Let $E=(fH)^{-1}(0)$ be a simple normal crossings divisor with
irreducible components $E_i$ (where $i\in\irr(E)$). 
The zero divisor of $fH$
can then be written as $\sum m_i E_i$, and the Jacobian ideal $\jac{H}$
of $H$ (which is the principal ideal of $\reg_Y$ 
characterized by $H^*\Omega^d_X = \jac{H} \Omega^d_Y$) can be
expressed as $\sum (n_i - 1) E_i$, where $n_i > 0$.
For $I\subset \irr(E)$ we define $E_I$ as the intersection 
$\bigcap_{i\in I} E_i$ and $E_I^\circ$ as $E_I - \bigcup_{j\not\in I}
E_j$
(for $I=\emptyset$ we get $E_\emptyset = Y$ and $E_\emptyset^\circ = Y- E$).
\begin{prop}
\label{formelzetanaiv}
We have the identity
\begin{gather*}\cS(f)(T) = \sum_{I\subseteq \irr(E)}
[E_I^\circ]_X\prod_{i\in I} \frac{\tate - 1}{T^{-m_i}\tate^{n_i} -1}.
\end{gather*}                 
in $\M_X[[T]]$.
\end{prop}
In particular, $\cS(f)$ lies in
$\M_X[(T^{-N}\tate^n-1)^{-1} \,|\,n,N\in\nat_{>0}]$.

The motivic zeta function also carries information about the
monodromy. 

Again, let $X$ be a smooth connected variety of dimension $d$ over
$k$, \linebreak and 
let $f:X\longrightarrow\aff{1}$
be a function on $X$. We will always denote the zero set of
$f$ by $X_0$.
For a natural number $n\ge 1$ we define

$$ S_n(f):=
\tate^{-nd}[\{\gamma\in\bog_n(X)\,|\, f\gamma(t)=t^n\}]_{X_0}
\in \M_{X_0}^{\muedach}.$$

Here the $\muedach$-action is induced by the $\mu_n$-action on $\De_n$
given by 
$t\mapsto \zeta t$.

Following Denef and Loeser, we define the 
\emph{motivic zeta function of $f$} as 
$$S(f)(T):=\sum_{n\ge 1}S_n(f)T^n \in \M_{X_0}^{\muedach}[[T]].$$

There is also a formula for $S(f)(T)$ in terms of an embedded
resolution of $X_0$ (for non-zero $f$):
Let $Y$, $n_i$, and $m_i$ be like above. For $I\subseteq \irr(E)$ we
denote the greatest common divisor $(m_i)_{i\in I}$ by $m_I$. 
For a point in $E_I^\circ$ in a neighborhood $U$ the function $fH$ 
can be written as $u\prod_{i\in I} x_i^{m_i}$, where $u$ is a unit and
$x_i$ is a local (analytic) coordinate defining $E_i$. We define 
$\widetilde{E_I^\circ}$ as the $\mu_{m_I}$-covering of $E_I^\circ$
given over $U\cap E_I^\circ$ by
$\{(z,p)\in\aff{1}\times U\cap E_I^\circ \sd z^{m_I} = u(p)^{-1}\}$.
These patch to a \linebreak $\mu_{m_I}$-covering of $E_I^\circ$:
If $y_i=\eta_ix_i$  
are other local (analytic) coordinates with $\eta_i$ units, and $f$ is
written as  $v\prod_{i\in I} y_i^{m_i}$, then $u=v\prod_{i\in
I}\eta_i^{m_i}$ and 
$$\{(z,p)\in\aff{1}\times U\cap E_I^\circ \sd z^{m_I} = u(p)^{-1}\}
\cong \{(z,p)\in\aff{1}\times U\cap E_I^\circ \sd z^{m_I} = v(p)^{-1}\}$$
via $(z,p)\mapsto (\prod_{i\in I} \eta_i^{\alpha_i}z,p)$, where
$\alpha_i:=\frac{m_i}{m_I}$.

The $\mu_{m_I}$-operation on $\widetilde{E_I^\circ}$ induces a good
$\muedach$-action over $E_I^\circ$.

There is also an intrinsic description of these coverings, compare
also~\cite{Looijenga}, Lemma 5.3 and~\cite{Lefschetz}, Lemma 2.5:

Denote by $\nu_{E_i}$ the normal bundle of $E_i$ in $Y$, denote the
complement of the zero section by $U_{E_i}$ and the fiber product of
the restrictions of the $(U_{E_i})_{i\in I}$ over $E_I^\circ$ by 
$U_I$. 
As $fH$ is a section of $\reg_Y(\sum_{i\in I} -m_i E_i)$, it induces
a morphism 
$$\bigotimes_{i\in I} \nu_{E_i}^{\otimes m_i}|_{E_I^\circ}
\longrightarrow \aff{1}.$$ The composition of this map with 
the morphism
$\prod_{i\in I} \nu_{E_i}|_{E_I^\circ}
\longrightarrow
\bigotimes_{i\in I} \nu_{E_i}^{\otimes m_i}|_{E_I^\circ}$ which sends
$(v_i)$ to $\otimes v_i^{\otimes m_i}$ restricted to $U_I$ induces
$U_I\longrightarrow \multgr$.

Let us spell this out in the analytic coordinates used above:
An element of the fiber of $U_I$ over $p$ given by
$(\lambda_i\frac{\partial}{\partial x_i})_{i\in I}$ is mapped to
$u(p)\prod\lambda_i^{m_i}$ (where $\lambda_i\in\multgr$).

Define $U_I(1)$ as the preimage of $1$ under the morphism 
$U_I\longrightarrow \multgr$. Denote $\sum_{i\in I} m_i$ by $c_I$. Then
$\muedach$ acts on $U_I(1)$ via the scalar action of $\mu_{c_I}$. 

We get a $\muedach$-equivariant mapping 
$U_I(1)\longrightarrow \widetilde{E_I^\circ}$, where
in the local coordinates the
$z$-coordinate is given as 
\begin{equation}
\label{formelz}
\prod_{i\in I} \lambda_i^{\alpha_i}
\end{equation}
(note that this is well defined, as for different coordinates $y_i=\eta_i
x_i$ with $\eta_i$ a unit we have 
$\frac{\partial}{\partial x_i}=\eta_i\frac{\partial}{\partial y_i}$ on
$E_i$).

As the $\alpha_i$ are relatively prime, the vector $(\alpha_i)$ can be
completed by $|I|-1$ more vectors to a basis of $\ganz^I$, so we
conclude that $U_I(1)\longrightarrow\widetilde{E_I^\circ}$ is a torus
bundle with fiber $\multgr^{|I|-1}$, which in turn establishes 
$\widetilde{E_I^\circ}\longrightarrow E_I^\circ$ as the covering obtained from
$U_I(1)\longrightarrow E_I^\circ$ by passing (fiber-wise) to connected
components.

\begin{rem}
\label{einheitnull}
Note that for the zero function $f=0$ we get
$S(0)(T)=0\in\M_X^{\muedach}[[T]].$
If on the other hand $f$ is a unit,
$X_0=\emptyset$ and consequently $\M_{X_0}=0$.
\end{rem}

Denef and Loeser prove the
\begin{thm}
\label{formelzeta}
We have the identity
$$S(f)(T)=\sum_{\emptyset\ne I \subseteq \irr(E)} 
(\tate-1)^{|I|-1}[\widetilde{E_I^\circ}]_{X_0}
\prod_{i\in I} \frac{1}{T^{-m_i}\tate^{n_i} -1}$$
in $\M_{X_0}^{\muedach}[[T]]$.
\end{thm} 

In particular, $S(f)$ lies in
$\M_{X_0}^{\muedach}[(T^{-N}\tate^n-1)^{-1} \,|\,n,N\in\nat_{>0}]$ and can
hence by evaluated at $T=\infty$. 
Following Denef and Loeser, 
we call $\psi_f:=-S(f)(\infty)\in \M_{X_0}^{\muedach}$ the
\emph{motivic nearby fiber of $f$}.
In terms of an embedded resolution of $X_0$ it is given by

\begin{equation}
\psi_f= \sum_{\emptyset\ne I\subseteq \irr(E)}
(1-\tate)^{|I|-1}[\widetilde{E_I^\circ}]_{X_0}. \label{formelnabij}
\end{equation}

\begin{rem}
\label{aufl}
Note that if $\pi:Y\longrightarrow X$ is an isomorphism outside $X_0$, 
we have $\psi_f = (\pi_0)_! \psi_{f\pi}$, where $\pi_0$ denotes the morphism 
$Y_0\longrightarrow X_0$ between the zero loci of $f\pi$ and $f$.
\end{rem}

\section{Relative equivariant Grothendieck groups of varieties}
\label{relequ}
In this section we collect some properties of relative equivariant 
Grothendieck groups of varieties.

\begin{rem}
\label{proj}
Suppose that $V\longrightarrow S$ is a vector
bundle of rank $n+1$ which carries a linear $G$-action over the action
on $S$. Denote
$\Proj(V)\longrightarrow S$ by $\nu$.
Then the endomorphism
$\nu_!\nu^*$ of $\Knull^G(\Var{S})$ is multiplication with $[\proj{n}]$.
\end{rem}

The presentation for $\Knull^G(\Var{S})$ given in~\cite{euler}
in case of a trivial $G$-action on the base variety $S$ can be
generalized slightly.
\begin{lem}
Suppose we are given a good $G\times H$-action on a variety $S$, such
that $G$ acts trivially and $H$ acts freely.  Then the morphism
$$\Knullgrstr{G\times H}(\Var{S})\longrightarrow \KnullGstr(\Var{H\raus
S})$$ which maps $[X]_S$ to $[H\raus X]_{H\raus S}$ is a
$\KnullGstr(\Vark)$-linear
isomorphism.
Furthermore it induces a $\KnullG(\Vark)$-linear isomorphism
$$\Knull^{G\times H}(\Var{S})\longrightarrow 
\Knull^G(\Var{H\raus S}).$$
\end{lem}

\begin{proof}
The $\KnullGstr(\Vark)$-linearity is quite clear.
Suppose we are given a variety $X\longrightarrow S$ with a good
$G\times H$-action over the action on $S$.  Then the $G\times
H$-equivariant mapping $X\longrightarrow H\raus X\times_{H\raus S} S$ is
an isomorphism (as both are \'{e}tale of the same degree over
$H \raus X$). Hence pulling back along $S\longrightarrow H\raus S$ is an
inverse for $\Knullgrstr{G\times H}(\Var{S})\longrightarrow
\KnullGstr(\Var{H\raus S})$.

A vector bundle $V\longrightarrow X$ with a linear $G\times H$-action over
a good action on $X$ descends to the vector bundle $H\raus V$ 
over $H\raus X$.
And as $\Proj(H\raus V)\cong H\raus \Proj(V)$ in this case, 
the above morphism
induces an isomorphism
$\Knull^{G\times H}(\Var{S})\longrightarrow 
\Knull^G(\Var{H\raus S})$
which is obviously $\KnullG(\Vark)$-linear.
\end{proof}
\begin{rem}
Applied to the case that $G$ is the trivial group this yields
$\Knullgrstr{H}(\VarS)\cong\Knull^H(\VarS)$ for a free good $H$-action
on $S$.
\end{rem}
\begin{rem}
The morphism $\Knull^{G\times H}(\Var{S})\longrightarrow 
\Knull^G(\Var{H\raus S})$
induces an $\M^G_k$-linear morphism
$$\M^{G\times H}_S\longrightarrow \M^G_{H\raus S}$$ which we will 
sometimes denote by $A\mapsto \overline{A}$.
\end{rem}

For convenience we spell out the presentations for 
$\Knullgrstr{G\times H}(\Var{S})$ and $\Knull^{G\times H}(\Var{S})$
we obtain from the above lemma. 

\begin{cor}
\label{presrelequfrei}
Suppose that we are given a good $G\times H$-action on $S$ such that
$G$ acts trivially and $H$ acts freely on $S$.
The group $\Knullgrstr{G\times H}(\VarS)$
has a 
presentation
as the abelian group generated by the isomorphism classes
of $S$-varieties with good $G\times H$-action over $S$ which are smooth over
$k$  and proper over $S$ subject 
to the relations $[\emptyset]_S=0$ and $[\Bl{Y}{X}]_S -[E]_S = [X]_S
-[Y]_S$, 
where $X$ is smooth over $k$ and proper over $S$ and carries a good
$G\times H$-action over $S$, 
$Y\subset X$ is a closed smooth $G\times H$-invariant 
subvariety, $\Bl{Y}{X}$ is the blow-up of $X$ along $Y$ and $E$ is the
exceptional divisor of this blow-up.
Moreover, we get the same group if 
we restrict to
varieties which are projective over $S$. 
We can also restrict to varieties such that $G\times H$ acts transitively on
the connected components.
\end{cor}

\begin{cor}
\label{erzprglfrei}
The group $\Knull^{G\times H}(\VarS)$ is the free abelian group on smooth
varieties, projective (respectively, proper) over $S$ with good
$G\times H$-action over $S$ (transitive on the connected components),
modulo blow-up relations and the subgroup generated by expressions of
the form $[G\times H \circlearrowright \Proj(V)]_S -[\proj{n}\times
(G\times H\circlearrowright X)]_S$, where $X$ is a smooth variety,
projective over $S$, with a $G\times H$-action transitive on the
connected components, and $V\longrightarrow X$ is a vector bundle of
rank $n+1$ over $X$ with a linear action over the action on $X$.
\end{cor}

\begin{cor}
We can define the duality endomorphism $\D_S^{G\times H}$ on $\M_S^{G\times H}$
as in~\cite{euler} and get the same formulae as developed there. 
Dualizing commutes with restriction and induction and with the morphism induced
by dividing out the free $H$-action.
\end{cor}

\begin{rem}
We have not used that $k$ is algebraically closed here.
\end{rem}

Instead of finite groups, we will also consider $\muedach$ 
and finite products of $\muedach$ with itself and with finite groups:

For a natural number $l$, 
we define a \emph{good ${\muedach}^l$-action} to be an
action coming from a good $(\mu_n)^l$-action. For a finite group $G$, 
we can also consider $G\times{\muedach}^l$-actions, which we call 
\emph{good} if they come from a good $G\times (\mu_n)^l$-action.
We define
$\Knullgrstr{{\muedach}^l}(\VarS)$, 
$\Knullgrstr{G\times {\muedach}^l}(\VarS)$ and $\Knull^{{\muedach}^l}(\VarS)$, 
$\Knull^{G\times {\muedach}^l}(\VarS)$ as above.
As they are direct limits of groups of the kind considered above, we get 
analogous presentations. We only spell out one.

\begin{cor}
\label{erzprglfreimuedach}
Suppose $G$ and $H$ are finite groups, and $l\in\nat$.
Suppose that $S$ carries a good $G\times H\times {\muedach}^l$ action,
such that $H$ acts freely and $G$ and ${\muedach}^l$ acts
trivially.
Then the group $\Knull^{G\times H\times{\muedach}^l}(\VarS)$ is the free 
abelian group on smooth
varieties, projective (respectively, proper) over $S$ with (good)
$G\times H\times {\muedach}^l$-action over $S$ 
(transitive on the connected components),
modulo blow-up relations and the subgroup generated by expressions of
the form $[G\times H \times{\muedach}^l \circlearrowright \Proj(V)]_S 
-[\proj{n}\times
(G\times H\times{\muedach}^l \circlearrowright X)]_S$, where $X$ is a 
smooth variety,
projective over $S$, with a $G\times H\times{\muedach}^l$-action transitive on the
connected components, and $V\longrightarrow X$ is a vector bundle of
rank $n+1$ over $X$ with a linear action over the action on $X$.
\end{cor}

\begin{rem}
As above, dividing out by $H$ induces an
$\M_k^{G\times{\muedach}^l}$-linear morphism
$$\M^{G\times H\times{\muedach}^l}_S
\longrightarrow \M^{G\times{\muedach}^l}_{H\raus S}$$ which we will 
also denote by $A\mapsto \overline{A}$.

Furthermore we also get a duality endomorphism $\D_S^{G\times H\times
{\muedach}^l}$ of $\M^{G\times H\times{\muedach}^l}_S$ which satisfies
the same relations as developed above.
\end{rem}

\section{The relative dual of an affine simplicial toric variety}
The aim of this section is the following
\begin{lem} \label{dualtorisch}
Let $X$ be an affine toric variety associated to a 
simplicial cone, let $X\longrightarrow S$ be proper, let $G$ be
a finite group acting on $X$ over $S$ via the torus, where $S$ carries
the trivial $G$-action.
Then $\D_S[X]_S=\tate^{-\dim X}[X]_S\in \M_S^G$.
\end{lem}

\begin{cor}
If $X$ is a toric variety associated to a simplicial fan, we have
$\D_X[X]_X=\tate^{-\dim X}[X]_X\in \M_X$.
\end{cor}

\begin{proof}
As for an open cover $\{U_i\}$ of $X$ the map
$\M_X\longrightarrow \prod_i \M_{U_i}$ is injective and commutes
with dualizing, we may assume that $X$ is defined by a simplicial cone.
\end{proof}

\begin{rem}
Note that if $X$ is complete, in particular 
$\D_k[X]=\tate^{-\dim X}[X]\in\M_k^G$.
\end{rem}

We need Lemma~\ref{formelglobal} on triangulations of simplices.
It is probably well known and perhaps should be proven by means of
toric geometry as it is derived from the Dehn-Sommerville equations 
which are the combinatorial counterpart of Poincar\'{e} duality for
toric varieties defined by complete simplicial fans.
First an auxiliary
\begin{lem}
$\sum_{l=0}^{n} \binom{n+1}{n-l} (t-1)^l = t^n+\dots+1$.
\end{lem}
\begin{proof}
This follows from 
$$\sum_{l=0}^{n} \binom{n+1}{l+1} (t-1)^{l+1} = ((t-1)+1)^{n+1}-1=t^{n+1}-1.$$
\end{proof}

For convenience we recall the
\emph{Dehn-Sommerville equations} --- see for example
\cite{FultonToric}, page 126.

\begin{thm}[Dehn-Sommerville equations]
Suppose we are given a triangulation of an 
$(m-1)$-sphere with $f_i$ faces of dimension
$i$. Let $f_{-1}:=-1$.
For $0\le p\le m$ we set
$$h_p=\sum_{i=p}^{m} (-1)^{i-p}\binom{i}{p} f_{m-1-i}.$$
Then $$h_p=h_{m-p}\text{ for }0\le p \le m.$$
\end{thm}

Suppose we are given a linear triangulation $S$ of an
$n$-simplex $\Delta$ which refines the standard 
triangulation $T$.
For $\sigma\in S$ we define $\sigma_{\Delta}$ to be the 
smallest simplex in $T$ which contains $\sigma$. We denote the
dimension of a simplex $\sigma$ by $|\sigma|$
and define the \emph{star of $\tau$ in $S$} by 
$\Star^S(\tau)=\{\sigma \in S \sd \tau\subseteq\sigma\}$.

\begin{lem}
\label{formellokal}
For a fixed $\tau\in S$ consider the polynomial 
$$g^S_{\tau}(t) = \sum_{\sigma \in \Star^S(\tau)} 
(-1)^{|\sigma_{\Delta}|} (t-1)^{|\sigma_{\Delta}|-|\sigma|}.$$

Then $g^S_{\tau}(t^{-1})=t^{|\tau|-n}g^S_{\tau}(t)$, 
in other words, $g^S_{\tau}$ is a polynomial of degree $\le n-|\tau|$ 
with symmetric coefficients. 
\end{lem}

\begin{proof}
We proceed by induction on $n-|\tau_{\Delta}|$.

Suppose $|\tau_{\Delta}|=n$. Then also $|\sigma_{\Delta}|=n$ for
all $\sigma \in \Star^S(\tau)$, hence 
\begin{align*}
(-1)^n g^S_{\tau}(t)&= \sum_{\sigma \in \Star^S(\tau)} 
(t-1)^{n-|\sigma|}\\ &= \sum_{k=|\tau|}^n c_k (t-1)^{n-k}
=\sum_{i=0}^{n-|\tau|} c_{n-i} (t-1)^i,
\end{align*}
where $c_k$ denotes the number of $k$-simplices in $\Star^S(\tau)$.
Without loss of generality we may assume that that $\Delta\subset
\R^n$ and that
$0\in \tau$. Denote the subspace generated by $\tau$ by $W_{\tau}$.
The projection of $\Star^S(\tau)$ to $\R^n/W_{\tau}$ generates a
complete fan whose intersection with a sphere around the origin gives
a triangulation with
$f_i=c_{i+|\tau|+1}$ simplices of dimension $i$. If we set $f_{-1}=1$
we have $f_i=c_{i+|\tau|+1}$ for all $-1\le i\le n-|\tau|-1$.
Hence 
\begin{align*}
(-1)^n g^S_{\tau}(t)&=\sum_{i=0}^{n-|\tau|} f_{n-|\tau|-1-i} (t-1)^i\\
&=\sum_{i=0}^{n-|\tau|}\sum_{p=0}^i 
f_{n-|\tau|-1-i} \binom{i}{p}t^p(-1)^{i-p}\\
&=\sum_{p=0}^{n-|\tau|}t^p\sum_{i=p}^{n-|\tau|} 
(-1)^{i-p}\binom{i}{p}f_{n-|\tau|-1-i}.
\end{align*}
The Dehn-Sommerville equations then imply that 
$g^S_{\tau}(t^{-1})= t^{|\tau|-n}g^S_{\tau}(t)$.

Now suppose $k:=|\tau_{\Delta}|<n$. Then $\tau$ is 
contained in a $k$-dimensional face of $\Delta$.
Without loss of generality we can assume that $\Delta\subset \R^n$
is the simplex spanned by $0$ and the $n$ standard basis vectors 
$e_1,\dots,e_n$, and that $\tau$ is contained in the facet spanned by
$0$ and $e_1,\dots,e_k$. Denote by $H$ the hyperplane spanned by
$e_1,\dots,e_{n-1}$. Denote the reflection at $H$ by $\rho$.
We get a new simplex $\Delta'=\Delta\cup \rho(\Delta)$ (spanned by
$e_1,\dots,e_{n-1}, e_n,-e_n$) with the linear triangulation
$S'=S\cup\rho(S)$ refining the standard triangulation of $\Delta'$
such that $|\tau_{\Delta'}|=|\tau_{\Delta}|+1$.

The star of $\tau$ in $S'$ decomposes as 
\begin{align*}
\Star^{S'}(\tau)&= \{\sigma\in \Star^S(\tau)\sd \sigma\subseteq H\}\\
             &\sqcup \{\sigma\in \Star^S(\tau)\sd \sigma \not\subseteq H\}
             \sqcup \rho\{\sigma\in \Star^S(\tau)\sd 
              \sigma \not\subseteq H\}.
\end{align*}
For $\sigma\in \Star^S(\tau)$ such that $\sigma\subseteq H$
we get $|\sigma_{\Delta'}|=|\sigma_{\Delta}|+1$ and 
$|\sigma_{H\cap \Delta}|=|\sigma_{\Delta}|$,
while for $\sigma\in \Star^S(\tau)$ such that 
$\sigma\not\subseteq H$ we have
$|\sigma_{\Delta'}|=|\sigma_{\Delta}|$.
Furthermore $|(\rho\sigma)_{\Delta'}|=|\sigma_{\Delta'}|$, and hence
\begin{align*}
g^{S'}_{\tau}(t)&=\sum_{\substack{\sigma\in \Star^S(\tau)\\
                \sigma\subseteq H}}
            (-1)^{|\sigma_{\Delta}|+1} 
            (t-1)^{|\sigma_{\Delta}|-|\sigma|+1}\\
          &+2 \sum_{\substack{\sigma\in \Star^S(\tau)\\ 
                \sigma \not\subseteq H}}
            (-1)^{|\sigma_{\Delta}|} 
            (t-1)^{|\sigma_{\Delta}|-|\sigma|}\\
&=2 g^S_{\tau}(t) - 2g^{H\cap S}_{\tau}(t) 
+ (1-t)g^{H\cap S}_{\tau}(t)\\
&=2 g^S_{\tau}(t) - (1+t) g^{H\cap S}_{\tau}(t).
\end{align*}
By the induction hypothesis we have 
$g^{S'}_{\tau}(t^{-1})= t^{|\tau|-n} g^{S'}_{\tau}(t)$ and 
$g^{H\cap S}_{\tau}(t^{-1})= t^{|\tau|-n+1} g^{H\cap S}_{\tau}(t)$.
Using the above equation we conclude
$g^S_{\tau}(t^{-1})= t^{|\tau|-n} g^S_{\tau}(t)$.
\end{proof}

\begin{lem}
\label{formelglobal}
Consider the polynomial 
$$h^S(t):= \sum_{\sigma\in S} (-1)^{|\sigma_{\Delta}|}
(t-1)^{|\sigma_{\Delta}|-|\sigma|} - 1.$$
Then $h^S(t^{-1})=t^{-(n+1)}h^S(t)$. 
\end{lem}

\begin{proof}
The cone $R$ on $S$ is a linear triangulation of 
an $(n+1)$-simplex. If $\tau$ is the top, we have
$g^{R}_{\tau}(t)= - h^S(t)$.
\end{proof}

\begin{proof}[Proof of Proposition~\ref{dualtorisch}]
Denote the fan defining $X$ by $\Tau$. 
We proceed by induction on $\dim \Tau$.

For $\dim \Tau = 0$ the claim holds, because in this case $X$ is smooth.

Now suppose $\dim \Tau \ge 1$.
Note that all toric constructions will be compatible with the action
of $G$. 
In the sequel we will denote the dimension of a cone $\tau$ by
$|\tau|$.
We have the orbit stratification
$X =\bigsqcup_{\tau\in \Tau} O_{\tau}$, where
$O_{\tau}$ is $(k-|\tau|)$-dimensional.

Furthermore, if $V_{\tau}$ denotes the closure of $O_{\tau}$, we have
the equation 
$[O_{\tau}]_S=\sum_{\tau'\supseteq\tau} 
(-1)^{|\tau'|-|\tau|}[V_{\tau'}]_S$ in
$\Knull^{G}(\VarS)$.

We choose a toric resolution of singularities $Y\longrightarrow X$.
It is given by a certain simplicial refinement $\Sigma$ of $\Tau$ and carries
the orbit stratification
$Y =\bigsqcup_{\sigma\in\Sigma} O_{\sigma}$.

For $\sigma\in\Sigma$ we denote by $\varphi(\sigma)\in \Tau$ the
smallest facet of $\Tau$ which contains $\sigma$.

Then $O_{\sigma}\longrightarrow O_{\varphi(\sigma)}$ is a (trivial)
$\multgr^{(|\varphi(\sigma)|-|\sigma|)}$-bundle.
As $G$  acts via the torus, 
$[O_{\sigma}]_S= (\tate-1)^{(|\varphi(\sigma)|-|\sigma|)} 
[O_{\varphi(\sigma)}]_S$ in $\Knull^{G}(\VarS)$.
Thus
\begin{align*}
[Y]_S &= \sum_{\tau\in \Tau}\sum_{\varphi(\sigma)=\tau} 
(\tate-1)^{|\tau|-|\sigma|} [O_\tau]_S\\
&= \sum_{\tau'\in\Tau} \biggl(\sum_{\tau\subseteq \tau'}
\sum_{\varphi(\sigma)=\tau} (-1)^{|\tau|}(\tate-1)^{|\tau|-|\sigma|}\biggr)
(-1)^{|\tau'|}[V_{\tau'}]_S\\
&=\sum_{\substack{\tau'\in\Tau\\|\tau'|\ge 1}}
 \biggl( \sum_{\tau\subseteq \tau'}
 \sum_{\varphi(\sigma)=\tau} (-1)^{|\tau|}(\tate-1)^{|\tau|-|\sigma|}\biggr)
 (-1)^{|\tau'|} [V_{\tau'}]_S + [X]_S
\end{align*}
For $\tau'\in\Tau$ such that $|\tau'|\ge 1$ consider
\begin{align*}
p^{\tau'}(t):&= \sum_{\tau\subseteq \tau'} \sum_{\varphi(\sigma)=\tau}
(-1)^{|\tau|}(t-1)^{|\tau|-|\sigma|}\\ 
&= \sum_{\substack{\tau\subseteq \tau'\\|\tau|\ge 1}} \sum_{\varphi(\sigma)=\tau}
(-1)^{|\tau|}(t-1)^{|\tau|-|\sigma|} + 1
\end{align*}
Now from Lemma~\ref{formelglobal} (applied to the intersection of
$\tau'$ with a transversal hyperplane and the triangulations
induced by $\Sigma$ and $\Tau$) we know that 
$p^{\tau'}(t^{-1})= t^{-|\tau'|} p^{\tau'}(t)$.
Furthermore $V_{\tau'}$ is an affine toric variety defined by a
simplicial cone of dimension $\dim \Tau -|\tau'|$, 
hence for $|\tau'|\ge 1$ we deduce from the induction 
hypothesis 
that 
$\D_S([V_{\tau'}]_S) 
= \tate^{-(\dim X -|\tau'|)}[V_{\tau'}]_S$.
Thus
$$\D_S([Y]_S - [X]_S) 
= \tate^{-\dim X}([Y]_S - [X]_S),$$
which completes the induction step.

\end{proof}

\section{A completion}
Let $f:X\longrightarrow \aff{1}$ be a non-zero function on a smooth
connected variety $X$, let $Y$ be an embedded resolution of
$X_0=f^{-1}(0)$ and let $E_I$ and $\widetilde{E_I^\circ}$ be as in
Section~\ref{motint}.  

\begin{defi}
Let $\emptyset\ne I \subset \irr(E)$. Define
$\widetilde{E_I}$ as the normalization of $E_I$ in 
$\widetilde{E_I^{\circ}}$.
\end{defi}

\begin{lem}\label{completion}
We have $\widetilde{E_I}|_{E_J}\cong\widetilde{E_J}$ for $I\subseteq J$
and $\D_{E_I}[\widetilde{E_I}]_{E_I}=\tate^{|I|-\dim
X}[\widetilde{E_I}]_{E_I}$ in $\M_{E_I}^{\muedach}$ (and in particular
$\D_{X_0}[\widetilde{E_I}]_{X_0}= \tate^{|I|-\dim
X}[\widetilde{E_I}]_{X_0}$ in $\M_{X_0}^{\muedach}$).
\end{lem}

\begin{proof}
Both statements are local on the embedded resolution of $X_0$, 
hence we can assume that  $f=ux_1^{m_1}\cdots x_k^{m_k}$, with $k\ge
1$ and $m_i>0$, where $x_1,\dots,x_k,\dots,x_n$ are local analytic
coordinates, and $I=\{1,\dots,l\}\subseteq J=\{1,\dots,l'\}$, 
where $1\le l\le l'\le k$.

Adjoining an $m_1$-th root of $u$ we again get a $\mu_{m_1}$-cover
$$\pi:Y:=\{(p,t)\in X\times \aff{1}\sd t^{m_1}=u(p)\}\longrightarrow X$$ 
with analytic coordinates 
$y_1=tx_1, y_2=x_2, \dots , y_n=x_n$ around $\{x_1=0\}$ such that 
$f\pi = y_1^{m_1}\dotsm y_k^{m_k}$.
Here $\xi\in\mu_{m_1}$  acts on $y_1$ by multiplication with $\xi$ and
trivially on $y_j$ for $j\ge 2$. Shrinking $X$ (outside $\{x_1=0\}$),
we may assume that $\nu=(y_1,\dots,y_n):Y\longrightarrow \aff{n}$ is
\'etale. Denote $\{y_1=\dots=y_l=0\}$ by
$F_I$ and $F_I-\bigcup_{j=l+1}^k\{y_j=0\}$ by $F_I^\circ$ (these are the
pullbacks of $E_I$ and $E_I^\circ$ to $Y$). Let 
$$\widetilde{F_I^\circ}:=\{(s,q)\in\aff{1}\times F_I^\circ\sd 
s^{m_I}=\prod_{j=l+1}^k y_j(q)^{-m_j}\}$$. 
Recall that 
$$\widetilde{E_I^\circ}=\{(z,p)\in\aff{1}\times E_I^\circ\sd
z^{m_I}=u(p)\prod_{j=l+1}^k x_j(p)^{-m_j}\}$$ and note that the
pullback of $\widetilde{E_I^\circ}$ to $F_I^\circ$ is isomorphic to
$\widetilde{F_I^\circ}$ via $(z,p,t)\mapsto (t^{\frac{m_1}{m_I}}z,
(p,t))$ (and that the same holds for $\widetilde{F^\circ_J}$).  Under
this isomorphism the $\muedach$-action comes from a
$\mu_{m_I}$-action, where $\zeta\in \mu_{m_I}$ acts on $s$ by
multiplication with $\zeta$, while $\xi\in\mu_{m_1}$ acts on $s$ by
multiplication with $\xi^{\frac{m_1}{m_I}}$.

Denote the normalization of $F_I$ in $\widetilde{F_I^\circ}$ by
$\widetilde{F_I}$. Consider the following cartesian diagram:
$$\xymatrix{\pi^*\widetilde{E_I}\ar[r]\ar[d]& \widetilde{E_I}\ar[d]\\
            F_I\ar[r]^\pi& E_I}$$
As $\pi$ is smooth and $\widetilde{E_I}$ is normal, 
$\pi^*\widetilde{E_I}$ is normal. Furthermore, it is finite and
surjective over $F_I$
and hence isomorphic to $\widetilde{F_I}$.
We conclude that $\widetilde{E_I}\cong \mu_{m_1}\raus\widetilde{F_I}$
and similarly for $\widetilde{E_J}$.
Hence it suffices that $\widetilde{F_I}|_{F_J}\cong\widetilde{F_J}$ for 
$I\subseteq J$ and that $\D_{F_I}[\widetilde{F_I}]_{F_I}=\tate^{|I|-\dim
X}[\widetilde{F_I}]_{F_I}\in\M_{F_I}^{\mu_{m_1}\times\muedach}$.

Now consider the \'etale morphism 
$\nu:Y\longrightarrow \aff{n}$. 
We denote the coordinates of $\aff{n}$ by
$w_1,\dots,w_n$, $\{w_1=\dots=w_l=0\}$ by $D_I$,  
$D_I-\bigcup_{j=l+1}^k\{w_j=0\}$ by $D_I^\circ$. Let 
$$\widetilde{D_I^\circ}:=\{(r,w)\in\aff{1}\times D_I^\circ\sd 
r^{m_I}=\prod_{j=l+1}^k w_j^{-m_j}\}$$
and denote the normalization of $D_I$ in $\widetilde{D_I^\circ}$ by 
$\widetilde{D_I}$.
Then we
have again a cartesian diagram:
$$\xymatrix{\widetilde{F_I}\ar[r]\ar[d]& \widetilde{D_I}\ar[d]\\
            F_I\ar[r]^{\nu}& D_I}$$
Thus it suffices that
$\widetilde{D_I}|_{D_J}\cong\widetilde{D_J}$ for $I\subseteq J$
and $\D_{D_I}[\widetilde{D_I}]_{D_I}=\tate^{|I|-\dim
X}[\widetilde{D_I}]_{D_I}\in\M_{D_I}^{\mu_{m_1}\times\muedach}$.
As the projection $p:\aff{n}\longrightarrow \aff{k}$ is smooth 
and $\D p^* = \tate^{k-n}p^*\D$, 
we may assume without loss of generality that $k=n$. In this case the claim
follows from Lemma~\ref{einschr} and Lemma~\ref{allgemein}.
\end{proof}

\begin{lem}
\label{einschr}
The restriction of the normalization $\widetilde{S}$ of 
$S:=\{s^d=x_1^{p_1}\dotsm x_k^{p_k}\}\subset \aff{1}\times\aff{k}$
to $\{x_1=0\} \subset \aff{k}$ is isomorphic to the normalization
$\widetilde{S'}$
of $S':=\{{s'}^{d'}=x_2^{p_2}\dotsm x_k^{p_k}\}\subset
\aff{1}\times\aff{k-1}$,
where $d'=(d,p_1)$. The $\mu_d$-action on $\widetilde{S}$ 
which is given by $\zeta(s,x_1,\dots,x_k)=(\zeta^{-1} s,x_1,\dots,x_k)$
restricts to the action induced by
$\zeta(s',x_2,\dots,x_k)=(\zeta^{-\frac{d}{d'}}s',x_2,\dots,x_k)$. 
In other words, it is given
via the canonical surjection
$\mu_d\longrightarrow\mu_{d'}$.  
\end{lem}

\begin{proof}
Let us first assume that the greatest common divisor
$(p_1,\dots,p_k,d)$ equals $1$ or equivalently that $S$ is
irreducible. 

Let $M$ be the lattice in $\R^k$ spanned by 
$\ganz^k$ and by $v:=(\frac{p_1}{d},\dots,\frac{p_k}{d})$. Then 
$\widetilde{S}\cong \Spec k[M^+]$, where $M^+=M\cap (\Rp)^k$.
The restriction to $\{x_1=0\}$ is then given as
$\Spec k[M^+_1]$, where $M_1:=\{\alpha \in M | \alpha_1 =
0\}$ and $M^+_1=M_1 \cap (\Rp)^k$ (compare~\cite{Toroidal}, page 16).

But $M'$ is generated by $e_2,\dots,e_k$ and $v':=(0,\frac{p_2}{d'},\dots,
\frac{p_k}{d'})$:

First note that $v'=\frac{d}{d'}v- \frac{p_1}{d'}e_1 \in M'$.

On the other hand, if $(\lambda_1,\dots,\lambda_k)+\mu v \in M'$, where
$\lambda_i, \mu \in \ganz$,
then $\lambda_1 + \mu\frac{p_1}{d} = 0$.  

Write $p_1=fd'$ and $d=ed'$ such that $(e,f)=1$. As $\mu\frac{f}{e}\in
\ganz$ also $\mu'=\frac{\mu}{e} \in \ganz$ and hence
$(\lambda_1,\dots,\lambda_k)+\mu v = (0,\lambda_2,\dots,\lambda_k)+
\mu'v'$ lies in the $\ganz$-module generated by $e_2,\dots, e_k$ and
$v'$.

Hence indeed $\widetilde{S}|_{\{x_1=0\}} \cong \widetilde{S'}$.  Note
that $s' = s^{\frac{d}{d'}}x_1^{-\frac{p_1}{d'}}$, hence the
$\mu_d$-action on $\widetilde{S'}$ comes from the $\mu_{d'}$-action
induced by $s'\mapsto {\zeta'}^{-1}s'$
via the homomorphism $\mu_d\longrightarrow\mu_{d'}$
which maps $\zeta$ to $\zeta^{\frac{d}{d'}}$, and the
$\mu_{p_1}$-action comes from the $\mu_{d'}$-action
via $\xi \mapsto \xi^{\frac{p_1}{d'}}$.

Now let us drop the assumption that 
$c=(p_1,\dots,p_k,d)$ equals $1$.
Note that $c=(d',p_2,\dots,p_k)$.
Let $e=\frac{d}{c}$, $q_i=\frac{p_i}{c}$, $e'=(e,q_1)$.
Note that $e'=\frac{d'}{c}$.
Let $\widetilde{T}$ be the normalization of 
$T=\{t^e=x_1^{q_1}\dotsm x_k^{q_k}\}$ and $\widetilde{T'}$ the
normalization of  
$T'=\{{t'}^{e'}=x_2^{q_2}\dotsm x_k^{q_k}\}$. They both carry a
$\mu_e$-action as above. 
The mapping
$$\mu_d\times^{\mu_e} T\longrightarrow S, 
\quad (\eta,t,x)\mapsto (\eta^{-1}t,x)$$
induces an isomorphism 
$$\mu_d\times^{\mu_e} \widetilde{T}\cong \widetilde{S}$$
over $\aff{k}$.
The $\mu_d$-action on $\widetilde{S}$ corresponds to the action 
on $\mu_d\times^{\mu_e} \widetilde{T}$ given by (left) multiplication 
on $\mu_d$.

Furthermore the mapping
$$\mu_d\times^{\mu_e} T'\longrightarrow S', 
\quad (\eta,t',x)\mapsto (\eta^{-\frac{d}{d'}}t',x)$$
induces an isomorphism 
$$\mu_d\times^{\mu_e} \widetilde{T'}\cong \widetilde{S'}.$$
As 
$$(\mu_d\times^{\mu_e} \widetilde{T})|_{\{x_1=0\}}=
\mu_d\times^{\mu_e} (\widetilde{T}|_{\{x_1=0\}})= 
\mu_d\times^{\mu_e} \widetilde{T'}$$ 
this establishes $\widetilde{S}|_{\{x_1 = 0\}} \cong \widetilde{S'}$. 

Now let us calculate the action of $\mu_d$ on 
$\mu_d\times^{\mu_e} \widetilde{T'}$ obtained from the action on 
$\mu_d\times^{\mu_e} \widetilde{T}$ and the corresponding action on 
$\widetilde{S'}$.

The action of $\mu_d$ is given by (left) multiplication on the first 
factor. It corresponds to the $\mu_d$-action on $\widetilde{S'}$
induced by $\zeta(s',x_2,\dots,x_k)=
(\zeta^{-\frac{d}{d'}}s',x_2,\dots,x_k)$. 
\end{proof}

\begin{lem}
\label{allgemein}
Let $p_1,\dots,p_k,d$ be natural numbers,
not necessarily relatively prime. Again denote by
$\widetilde{S}$ the normalization of 
$S:=\{t^d=x_1^{p_1}\dotsm x_k^{p_k}\}\subset \aff{1}\times\aff{k}$.
Then $$\D_{\aff{k}}([\widetilde{S}]_{\aff{k}}) 
= \tate^{-k}[\widetilde{S}]_{\aff{k}}\text{ in
}\M_{\aff{k}}^{\mu_d\times \mu_{p_1}}.$$
Here the $\mu_d$-action is 
induced by $\zeta(t,x_1,\dots,x_k)=(\zeta^{-1} t,x_1,\dots,x_k)$, and
the \linebreak 
$\mu_q$-action is given by the $\mu_d$-action on $\widetilde{S}$ via
$\mu_q\longrightarrow \mu_d$, $\xi \mapsto \xi^{\frac{q}{d}}$.
\end{lem}

\begin{proof}

As dualization commutes with restriction, it is enough to prove the claim
in $\M_{\aff{k}}^{\mu_d}$. 

If the greatest common divisor $(p_1,\dots,p_k,d)$ equals $1$, the
claim follows directly from Lemma~\ref{dualtorisch}. 

The general case follows from the fact that dualizing commutes with
induction as in the proof of Lemma~\ref{einschr}.
\end{proof}

\section{The dual of the motivic nearby fiber and two functional equations}
Let $X$ be a smooth connected variety of dimension $d$, let
$f:X\longrightarrow \aff{1}$ be a (non-zero) function.
\begin{thm}\label{selfdual}
We have
$\D_{X_0} \psi_f = \tate^{1-d} \psi_f$ in $\M_{X_0}^{\muedach}$.
\end{thm}

So $\psi_f$ behaves like a smooth $(d-1)$-dimensional variety, proper
over $X_0$.

\begin{rem}\label{zerlproj}
Let $J$ be a finite nonempty set. Then
$$\sum_{\emptyset \ne I \subseteq J} (\tate-1)^{|I|-1} 
= [\proj{|J|-1}].$$
\end{rem}

\begin{proof}[Proof of Theorem~\ref{selfdual}]

According to Formula~(\ref{formelnabij}) we have
$$\psi_f= \sum_{\emptyset\ne I \subseteq \irr(E)}
(1-\tate)^{|I|-1}[\widetilde{E^\circ_I}]_{X_0}.$$ 
Due to Lemma~\ref{completion} $[\widetilde{E^\circ_I}]_{X_0} = \sum_{J\supseteq I} (-1)^{|J|-|I|} 
[\widetilde{E_J}]_{X_0}$, and hence 
\begin{align*}
\psi_f&= \sum_{\emptyset\ne I \subseteq \irr(E)} 
   (1-\tate)^{|I|-1}\sum_{J\supseteq I} (-1)^{|J|-|I|}[\widetilde{E_J}]_{X_0}\\
&=\sum_{\emptyset\ne J \subseteq\irr(E)}
(-1)^{|J|-1}[\widetilde{E_J}]_{X_0}
\sum_{\emptyset\ne I \subseteq J} (\tate-1)^{|I|-1}\\
&=\sum_{\emptyset\ne J \subseteq\irr(E)}
(-1)^{|J|-1}[\proj{|J|-1}][\widetilde{E_J}]_{X_0}\text{ due to
Remark~\ref{zerlproj}.}
\end{align*} 
Again using Lemma~\ref{completion} we conclude that
$$\D_{X_0}(\psi_f)=  
\sum_{\emptyset\ne J \subseteq\irr(E)}
(-1)^{|J|-1}\tate^{-|J|+1}[\proj{|J|-1}]\tate^{-d+|J|}[\widetilde{E_J}]_{X_0}
= \tate^{-d+1}\psi_f.$$
\end{proof}
Recall that Proposition~\ref{formelzetanaiv} shows that
$\cS(f)$ lies in
$\M_X[(T^{-m}\tate^n-1)^{-1} \,|\,n,m\in\nat_{>0}]$. 

Consider the ring 
$\Pl_k:=\M_k[T,T^{-1},(T^{-m}\tate^n-1)^{-1} \,|\,n,m\in\nat_{>0}]$.
The duality involution on $\M_k$ can be extended to a ring involution
$\D^{\Pl}_k$ 
of $\Pl_k$ by setting $\D^{\Pl}_k(T)=T^{-1}$.

We set
$\Pl_X:=\M_X\otimes_{\M_k} \Pl_k =
\M_X[T,T^{-1},(T^{-m}\tate^n-1)^{-1} \sd n,m\in\nat_{>0}]$.
The duality involution $\D_X$
induces a $\D^{\Pl}_k$-linear involution $\D^{\Pl}_X$.

\begin{claim}\label{funktnaiv}
The functional equation 
$\D^{\Pl}_X \cS(f) = \tate^{-d}\cS(f)$
holds in $\Pl_X$.
\end{claim}
\begin{proof}
This is much easier to see than the functional equation in 
\cite{MotivicIgusa}. The calculation given here can already been found
there.

Note that $[E^{\circ}_I]_X=\sum_{J\supseteq I} (-1)^{|J|-|I|} [E_J]_X$.

Furthermore for a finite collection of elements $(a_l)_{l\in L}$ in a
commutative ring the identity $\prod_{l\in L}(a_l -1) =
\sum_{K\subseteq L} (-1)^{|L|-|K|}\prod_{k\in K} a_k$ holds.

This together with Proposition~\ref{formelzetanaiv} yields

\begin{gather*}\cS(f) = \sum_{J\subseteq \irr(E)}
[E_J]_X\prod_{j\in J} B_j,
\end{gather*}                 
where $B_j = \frac{\tate -1}{T^{-m_j}\tate^{n_j}-1} -
1=\frac{\tate-T^{-m_j}\tate^{n_j}}{T^{-m_j}\tate^{n_j}-1}$.
Now note that $\D^{\Pl}_k(B_j) = \tate^{-1}B_j$ and that 
$\D_X[E_J]_X = \tate^{|J|-d}[E_J]_X$.
\end{proof}

Theorem~\ref{formelzeta} and Lemma~\ref{completion} yield
\begin{multline*}
S'(f)(T):= (\tate-1)S(f)(T) + \sum_{\emptyset\ne
J\subseteq\irr(E)} (-1)^{|J|}[\widetilde{E_J}]_{X_0}\\
=\sum_{\emptyset\ne J\subseteq
\irr(E)}[\widetilde{E_I}]_{X_0}\sum_{I\subseteq J} (-1)^{|J|-|I|}\prod_{i\in
I} \frac{\tate-1}{T^{-m_i}\tate^{n_i} -1}, 
\end{multline*}
and as in the proof of Claim~\ref{funktnaiv} we
conclude that 
$$\D_{X_0}^{\Pl} S'(f) = \tate^{-d} S'(f)$$
in 
$\M_{X_0}^{\muedach} [T,T^{-1},(T^{-m}\tate^n-1)^{-1}\sd n,m\in\nat_{>0}]$, where $\D^{\Pl}_{X_0}(T)=T^{-1}$.
\section{Some properties of motivic zeta functions 
and the motivic nearby fiber}
Let again $X$ be smooth connected variety over $k$ and $f:X\longrightarrow 
\aff{1}$ a function on $X$. We denote the zero locus of $f$ by $X_0$.

\begin{defi}
For a morphism $\pi:X'\longrightarrow X$ of varieties we will denote
the zero locus of $f\pi$ by $X'_0$ and the induced morphism
$X'_0\longrightarrow X_0$ by $\pi_0$.
\end{defi}

We collect some properties of the motivic zeta function and the motivic 
nearby fiber. For this purpose we first introduce a notation.
\begin{defi}
Let $m\ge 1$ be a natural number. Let $X$ be a variety with good
$\muedach$-action. Then we denote by
$\Ind^{(m)}:\M^{\muedach}_X\longrightarrow\M^{\muedach}_X$ 
the morphism induced by $\Ind^{\mu_{nm}}_{\mu_n}$, where
$\mu_n\longrightarrow \mu_{nm}$ is the inclusion $\zeta\mapsto \zeta$.
\end{defi}

\begin{listprop} \label{eigenschaftenzeta}
Suppose $u:X\longrightarrow \multgr$ is a morphism. Then
after a finite \'{e}tale base change $\pi:\widetilde{X}\longrightarrow X$ 
(taking a sufficiently high root of $u$) we get 
$S(f\circ\pi)(T)=S((uf)\circ\pi)(T)$ and in particular
\begin{equation}\label{einheit}
\psi_{f\circ\pi}=\psi_{(uf)\circ\pi}.
\end{equation}

Suppose $\pi:X'\longrightarrow X$ is a smooth morphism of smooth
connected varieties. Then
$S(f\pi)(T)={\pi_0}^*S(f)(T)$ and in particular
\begin{equation}
\label{glatt}
\psi_{f\pi}={\pi_0}^*\psi_f.
\end{equation}

For a natural number $m\ge 1$ we have $S(f^m)(T)=\Ind^{(m)} S(f)(T^m)$ and in
particular 
\begin{equation}
\label{ind}
\psi_{f^m}=\Ind^{(m)}\psi_f.
\end{equation}
\end{listprop}

\begin{proof}
The first equality follows directly from Theorem~\ref{formelzeta}.

For the second identity, note that
$\bog_n(X')\longrightarrow \pi^*\bog_n(X)$ is a locally trivial
fibration with fiber $\aff{nd}$, where $d=\dim X'-\dim X$.
Therefore $S_n(f\pi)={\pi_0}^*S_n(f)$.

Finally,
\begin{align*}
S_{nm}(f^m)&=\tate^{-nm\dim X}
[\gamma\in\bog_{nm}(X)\sd f^m\gamma(t)=t^{nm}]_{X_0}\\
&=\tate^{-nm\dim X}\tate^{n(m-1)\dim X}[\gamma\in\bog_{n}(X)\sd f\gamma(t)=\zeta t^n  \text{  where }\zeta\in\mu_m]_{X_0}\\
&=\Ind^{(m)}S_n(f).
\end{align*}
\end{proof}

\begin{rem}
Equality~(\ref{einheit}) does \emph{not} always hold before taking a base
change. For example, let $k=\C$.
Consider the open subvariety $U\subset\aff{1}$ which is 
the complement of the zero locus of $g(x)=x^3+ax+b$, where we assume that 
$g$ has no multiple roots. On $U\times \aff{1}$ consider the functions 
$f(x,y)=y^2$ and $u(x,y)=g(x)$. Then $\psi_f=\mu_2\times U\in\M_U^{\muedach}$ 
and $\psi_{uf}=\{t^2=u(x)\}\in\M_U^{\muedach}$ (set $t=z^{-1}$ in
Formula~\ref{formelnabij}). 
Their images in $\M_{\C}$ are
not equal: They are distinguished by the Hodge character.
\end{rem}

We will need the following lemma later on to define a nearby cycle morphism.

\begin{lem} \label{aufbl}
Let $X$ be a smooth connected variety and $Y\subset X$ a smooth
closed subvariety, let $f:X\longrightarrow \aff{1}$ be a function on
$X$. Let $\pi:\Bl{Y}{X} \longrightarrow X$ be the blow-up of $X$ along
$Y$, let $E$ be the exceptional divisor of the blow-up. Denote the
inclusion $Y\hookrightarrow X$ by $\iota$ and the inclusion $E\hookrightarrow
\Bl{Y}{X}$ by $\iota'$. Let $g:=f\iota$,
$f'=f\pi$ and 
$g'=f'\iota'$. Then
$$\psi_f - {\iota_0}_!\psi_g = 
{\pi_0}_!\psi_{f'}-{\pi_0}_!{\iota'_0}_!
\psi_{g'}.$$

\end{lem} 

\begin{proof}
Note that if $Y\subseteq f^{-1}(0)$ the claim follows from 
Remarks~\ref{einheitnull} and ~\ref{aufl}. Hence we may assume that
$Y$ is not contained in the zero locus of $f$.

If $f^{-1}(0)$ is a simple normal crossings divisor which has normal
crossings with $Y$, the same holds for ${f'}^{-1}(0)$ and $E$, and
Formula~\ref{formelnabij} yields
$$\psi_g={\iota_0}^*\psi_f\text{ and }
\psi_{g'}={\iota'_0}^*\psi_{f'}.$$
Hence $\psi_f-{\iota_0}_!\psi_g = {j_0}_!{j_0}^*\psi_f$, where $j$
denotes the inclusion $X-Y\hookrightarrow X$. Now ${j_0}^*\psi_f =
\psi_{fj}$, hence $\psi_f-{\iota_0}_!\psi_g = {j_0}_!\psi_{fj}$ and
similarly for $\psi_{f'}-{\iota'_0}_!\psi_{g'}$, thus
the claim follows from the fact that $\pi$ is an isomorphism outside
$Y$.

In the general case we first choose an embedded resolution $X^\natural$ of the zero
locus of $f$. Denote the closure of the inverse image of $Y-f^{-1}(0)$
by $Y^\natural$, denote the
function $X^\natural\longrightarrow X \longrightarrow \aff{1}$ by $f^\natural$.
Now we choose an embedded resolution $\widehat{Y}\subset\widehat{X}$ 
of $Y^\natural$ which is compatible with the zero divisor of $f^\natural$.

The situation is as follows:

$$
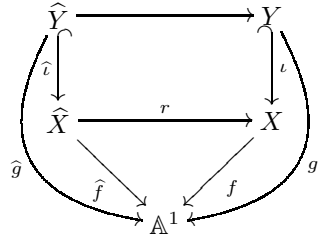
\xymatrix{
\widehat{Y}\ar[rr]\ar@{^{(}->}[d]_{\widehat{\iota}}
\ar@/_8ex/[ddr]_{\widehat{g}}
&& {Y}\ar@{_{(}->}[d]^{\iota}\ar@/^8ex/[ddl]^g\\
\widehat{X}\ar[rr]^r\ar[dr]_{\widehat{f}}&&X \ar[dl]^f\\
&\aff{1}&}$$

Here $\widehat{f}^{-1}(0)$ is a simple normal crossings divisor and has
normal crossings with $\widehat{Y}$.
Note that $r$ induces an isomorphism outside the zero locus of $f$,
hence in the diagram
$$\xymatrix{\Bl{\widehat{Y}}{\widehat{X}}\ar^{\widehat{\pi}}[d]\ar@{-->}[r]^{\varphi}
&\Bl{Y}{X}\ar[d]\\
\widehat{X}\ar[r]&X}$$
the proper birational map $\varphi$ induces an isomorphism outside the
zero locus of $f'$. Thus we can find $\widetilde{X}$ smooth making
$$\xymatrix{&\widetilde{X}\ar[dl]_p\ar[dr]^q&\\
\Bl{\widehat{Y}}{\widehat{X}}\ar@{-->}[rr]^{\varphi}&&\Bl{Y}{X}}$$
commutative such that the closure $\widetilde{E}$ of the inverse image
of $E-(f')^{-1}(0)$
is smooth and such that $p$ and $q$ are proper and isomorphisms
outside the zero loci of $f'$ and
$\widehat{f}'=\widehat{f}\widehat{\pi}$ 
(we can take a
resolution of singularities of the graph of $\varphi$ and then an
embedded resolution of the closure of the inverse image of 
$E-(f')^{-1}(0)$). We denote the inclusion  
$\widetilde{E}\hookrightarrow\widetilde{X}$ by $\widetilde{\iota}$,
$f'q$ by $\widetilde{f}$ and
$\widetilde{f}\widetilde{\iota}$ by $\widetilde{g}$. Furthermore
we define $\widehat{\iota'}$ as the inclusion $\widehat{E}\hookrightarrow
\Bl{\widehat{Y}}{\widehat{X}}$ and
$\widehat{g}':=\widehat{f}'\widehat{\iota}'$. 
From the commutative diagram
$$\xymatrix{
&\widetilde{X}\ar[dl]_p\ar[dr]^q\ar'[dd][ddd]_(-1){\widetilde{f}}&\\
\Bl{\widehat{Y}}{\widehat{X}}
\ar@/_8ex/[ddr]_{\widehat{f}'}\ar[d]^{\widehat{\pi}}&&
\Bl{Y}{X}\ar@/^8ex/[ddl]^{f'}\ar[d]_{\pi}\\
\widehat{X}\ar[dr]_{\widehat{f}}\ar[rr]^(.4){r}&& X\ar[dl]^f\\
&\aff{1}&}$$
we conclude:
\begin{alignat*}{2} 
\psi_f-\iota_{0!}\psi_g &= 
r_{0!}\psi_{\widehat{f}}-r_{0!}\widehat{\iota}_{0!}\psi_{\widehat{g}}
&\quad&\text{due to Remark \ref{aufl}}\\
&= r_{0!}\widehat{\pi}_{0!}\psi_{\widehat{f}'}-
r_{0!}\widehat{\pi}_{0!}\widehat{\iota}'_{0!} \psi_{\widehat{g}'} 
&\quad&\text{due to the above discussion}\\
&= r_{0!}\widehat{\pi}_{0!}p_{0!}\psi_{\widetilde{f}}-
r_{0!}\widehat{\pi}_{0!}p_{0!}\widetilde{\iota}_{0!}
\psi_{\widetilde{g}}
&\quad&\text{due to Remark \ref{aufl}}\\
&=\pi_{0!}q_{0!}\psi_{\widetilde{f}} -
\pi_{0!}q_{0!}\widetilde{\iota}_{0!}\psi_{\widetilde{g}}\\
&= \pi_{0!}\psi_{f'}
-\pi_{0!}{\iota'}_{0!} \psi_{g'}
&\quad&\text{due to Remark \ref{aufl}}
\end{alignat*}
\end{proof}

Now suppose there is a good action of a finite group $G$ on $X$
which is transitive on the connected components of $X$ and leaves $f$ 
invariant.

Hence we can regard $\cS(f)(T)$ as an element of 
$\M^G_X[[T]]$. The transformation formula also holds in the
equivariant setting (as an element of $g\in G$ induces an affine
action over the base in the fibrations of the Key lemma~9.2 
in~\cite{Looijenga} and on the $\aff{\dim X}$-bundle 
$\bog_{n+1}(X)\longrightarrow \bog_n(X)$). 
Hence we get a formula analogous to Proposition~\ref{formelzetanaiv} 
in the equivariant context if we choose an $G$-equivariant embedded
resolution of the zero locus (where the summation runs over the
\emph{orbits} of finite subsets of $\irr{E}$).

The $\muedach$-action on $\bog_n(X)$ induced by $\mu_n$ 
commutes with the action of
$G$ and hence $S(f)(T)$ can be regarded as an element of 
$\M^{G\times \muedach}_{X_0}[[T]]$ and $\psi_f$ as an element of 
$\M^{G\times \muedach}_{X_0}$. 

Also Theorem~\ref{formelzeta} and Formula~(\ref{formelnabij}) 
have analogues in this context, 
if we choose an equivariant embedded resolution of
the zero locus (from the intrinsic description of
$\widetilde{E_I^\circ}$ it is also clear that 
$\bigsqcup_{g\in G/\Stab_G(I)} \widetilde{E_{gI}^\circ}$ carries a $G$-action, where $\Stab_G(I)$
denotes the stabilizer of $I$ in $G$).

Lemma~\ref{aufbl} and properties~\ref{eigenschaftenzeta} also hold in the
equivariant setting.

We can replace $G$ by $G\times{\muedach}^l$ and get the same identities
as before.

\begin{rem}
\label{quotient}
Suppose we are given a good \emph{free} action of a finite group $H$ 
which is transitive on the
connected components on a smooth variety $X$. 
Suppose we have a $H$-invariant function $f$ on $X$. It then induces a
function $\overline{f}$ on $\overline{X}=H\raus X$. 
Furthermore $\bog_n(\overline{X})=H\raus \bog_n(X)$ and 
$S_n(\overline{f}) = H\raus S_n(f)$,
hence $S(\overline{f})(T) = \overline{S(f)(T)}$ (we extend the quotient 
morphism by $T\mapsto T$)
and in particular
$$\psi_{\overline{f}} = \overline{\psi_{f}}.$$
\end{rem}

\section{The nearby cycle morphism}
Let $X$ be a (not necessarily smooth) variety over $k$, 
let $f:X\longrightarrow \aff{1}$ be a function. 
Let $X_0:=f^{-1}(0)$. We want to define a
\emph{nearby cycle morphism}
$$\Psi_f: \M_X\longrightarrow \M_{X_0}^{\muedach}$$ such that
for a proper morphism
$\pi:X'\longrightarrow X$ we get ${\pi_0}_! \Psi_{f\pi}=\Psi_f \pi_!$
and furthermore in the case of a smooth connected variety $X$
the image of $\eins_X$ is $\psi_f$. 
For this purpose we define $\Psi_f$ on $\Knull(\Var{X})$ first.
\begin{defi}
Let $p:Y\rightarrow X$ be a proper morphism, where $Y$ is a
smooth connected $k$-variety.
Then we set $\Psi_f([Y]_X):= {p_0}_!(\psi_{fp})$.
\end{defi}
\begin{claim}
The morphism $\Psi_f$ is compatible with the blow-up relations (and
hence well defined).
\end{claim}
\begin{proof}
Let $p:Y\longrightarrow X$ as above, let $Z\subset Y$ be a smooth
connected closed subvariety. 
Denote by $\pi:Y'\longrightarrow Y$ the blow-up of $Y$ along
$Z$, denote the exceptional divisor by $E$.
Denote the inclusion
$Z\hookrightarrow Y$ by $\iota$ and the inclusion $E\hookrightarrow
Y'$ by $\iota'$. 
Then we have
\begin{align*}
\Psi_f([Y]_X)-\Psi_f([Z]_X)&={p_0}_!\psi_{fp}-(p_0\iota_0)_!\psi_{fp\iota}\\
&=(p_0\pi_0)_!\psi_{fp\pi}-(p_0\pi_0\iota'_0)_!\psi_{fp\pi\iota}
\text{ due to Lemma~\ref{aufbl}}\\
&=\Psi_f([Y']_X)-\Psi_f([E]_X).
\end{align*}
\end{proof}

\begin{claim}
The morphism $\Psi_f$ is $\Knull(\Vark)$-linear.
\end{claim}
\begin{proof}
Let $W$ be a smooth complete variety over $k$, 
let $p:Y\longrightarrow X$ be as above.
Denote the projection $W\times Y\longrightarrow Y$ by $\pi$. 
Then 
\begin{align*}
\Psi_f([W\times Y]_X)&=(p_0\pi_0)_!\psi_{fp\pi}
={p_0}_!{\pi_0}_!\pi_0^*\psi_{fp} \text{ due to Formula~(\ref{glatt}})\\
&={p_0}_!([W]\psi_{fp})=[W]\Psi_f([Y]_X).
\end{align*}
\end{proof}

Hence $\Psi_f$ can be extended to an $\M_k$-linear
morphism $$\Psi_f:\M_X\longrightarrow \M_{X_0}^{\muedach}.$$ 

\begin{listprop}
\label{eigenschaftenfunctor}
For $\pi:X'\longrightarrow X$ proper we get ${\pi_0}_! \Psi_{f\pi}=\Psi_f \pi_!$.

If $\pi:X'\longrightarrow X$ is a smooth morphism, 
then $\Psi_{f\pi}\pi^*={\pi_0}^*\Psi_f$.

If $\iota:X_0 \hookrightarrow X$ denotes the inclusion,
$\Psi_f \iota_!=0$.

For a natural number $m\ge 1$ we get $\Psi_{f^m}=\Ind^{(m)}\Psi_f$.

Dualizing and the nearby cycle morphism commute up to a factor: 
$\D^{\muedach}_{X_0} \Psi_f=\tate\Psi_f \D_X$.
\end{listprop}

\begin{proof}
The first identity holds by construction of $\Psi_f$.

For the second formula, suppose that
$\pi:X'\longrightarrow X$ is a smooth morphism.
Let $p:Y\longrightarrow X$ be a proper morphism, where $Y$ is a smooth
variety over $k$. 
Let 
$$\xymatrix{Y'\ar[r]^{p'}\ar[d]_{\pi'}&X'\ar[d]^{\pi}\\
            Y\ar[r]^p& X}$$
be cartesian. Then $\pi'$ is smooth, and $p'$ is proper, furthermore
$$\xymatrix{Y'_0\ar[r]^{p'_0}\ar[d]_{\pi'_0}&X'_0\ar[d]^{\pi_0}\\
            Y_0\ar[r]^{p_0}& X_0}$$
is cartesian, too.
We therefore get
\begin{align*}
\Psi_{f\pi}(\pi^*([Y]_X))&=
{p'_0}_!\psi_{fp\pi'}={p'_0}_!{\pi'_0}^*\psi_{fp}\text{ due to
Formula~(\ref{glatt})} \\ 
&= {\pi_0}^*{p_0}_!\psi_{fp} = {\pi_0}^*\Psi_f([Y]_X).
\end{align*}

The third identity follows from the fact that $S(0)(T)=0$ 
(see Remark~\ref{einheitnull}).

The last two identities follow from
Formula~(\ref{ind}) and Theorem~\ref{selfdual},
respectively. 
\end{proof}

We can also define $\Psi_f$ in the equivariant setting.
Let $G$ be a group of the form $H\times H'\times {\muedach}^l$, where
we assume that $H$ acts freely 
and $H'\times {\muedach}^l$ acts
trivially on the base variety $X$. We assume that we have a
$G$-invariant function $f:X\longrightarrow \aff{1}$.
We first get a $\KnullGstr(\Vark)$- linear morphism 
$$\Psi_f:\KnullGstr(\Var{X})\longrightarrow \M^{G\times\muedach}_{X_0}$$
as before.

\begin{claim}
Suppose $Y$ is a smooth variety, projective
over $X$ and that $V\longrightarrow Y$ is a vector
bundle of rank $n+1$ with a linear $G$-action over the 
action on $X$. Then 
$$\Psi_f([G\circlearrowright\Proj(V)]_X)
=\Psi_f([\proj{n}\times(G\circlearrowright Y)]_X).$$ 
\end{claim}

\begin{proof}
Let $\nu:\Proj(V)\longrightarrow Y$ be induced by the structure map
of $V$. Denote the morphism $Y\longrightarrow X$ by $p$. 
We have
\begin{align*}
\Psi_f([\Proj(V)]_X)
&= (p_0\nu_0)_!\psi_{fp\nu}
= (p_0\nu_0)_!(\nu_0)^*\psi_{fp}\text{ due to Formula~(\ref{glatt})}\\
&= {p_0}_!([\proj{n}]\psi_{fp})\text{ due to Remark~\ref{proj}}\\
&= [\proj{n}]\Psi_f([Y]_X) 
= \Psi_f([\proj{n}\times Y]_X).
\end{align*}
\end{proof}

Hence $\Psi_f$ induces a $\KnullG(\Vark)$-linear morphism
$$\Psi_f:\KnullG(\Var{X})\longrightarrow \M^{G\times\muedach}_{X_0}$$
and an $\M_k^G$-linear morphism
$$\Psi_f:\M^G_X\longrightarrow \M^{G\times\muedach}_{X_0}.$$ 

The first four identities of the List~\ref{eigenschaftenfunctor} also
hold in the equivariant setting. We do not know the relationship
between dualizing and the nearby cycle morphism in this context.

Furthermore, from Remark~\ref{quotient} we conclude
\begin{prop}
\label{quotientfunctor}
Suppose we have a good $G\times H$-action (transitive on the connected
components) on the smooth variety $X$ and a $G\times H$-invariant
function $f:X\longrightarrow \aff{1}$. 
Suppose that $H$ is a finite group and acts freely. Denote the
function induced on $\overline{X}= H\raus X$ by $\overline{f}$.
Then for $A\in\M_X^{G\times H}$ we have 
$$\overline{\Psi_f(A)}=\Psi_{\overline{f}}(\overline{A})
\in\M_{\overline{X}_0}^{G\times \muedach}.$$
\end{prop}

\end{document}